\documentclass[12pt]{article}
\usepackage{amssymb,amsmath}
\newtheorem{theorem}{Theorem}

\newtheorem{proposition}{Proposition}
\newtheorem{corollary}{Corollary}
\newtheorem{definition}{Definition}

\begin{document}

\title{
        Conservation laws for a class of Third order \\
        \noindent Evolutionary differential systems} 

\author{ Sung Ho Wang \\ }
\date{}
\maketitle

\noindent \textbf{\Large{Introduction} }

\noindent For the second order scalar evolutionary differential system (equation), it is known that every conservation law can be defined on the second jet space. [BG2]  \, In fact, there is an upper bound on the order of the conservation laws for the even order scalar evolution equation. However, it is not true in case of the odd order equations, as is illustrated by \emph{KdV} equation, and the infinite prolongation becomes necessary.

The lowest order conservation laws of the general quasilinear third order scalar evolution
equations are considered in [F], where the normal forms for the equations admitting such conservation laws are found. [K] considers the evolution equations with higher order conservation laws, and [M] provides an extensive list of the equations that possess some kind of formal symmetry.  

The purpose of this paper is to study the conservation laws of the third order scalar evolution equations of the form
\begin{equation} 
u_t = f( x, u, u_x ) u_{xxx} + g( x, u, u_x, u_{xx} )  
\end{equation}
via exterior differential systems and characteristic cohomology. The differential system and characteristic cohomology approach not only facilitates the computations involved, but it naturally gives rise to the convenient notion of the \emph{weight} of a conservation law, which will replace the order of a conservation law in our treatment. We find a universal integrability condition for an equation to have a higher order (weight) conservation law, and the computation of examples also suggests the following conjecture.

\emph{
If the differential system (equation) has two conservation laws of distinct weights $\geq \, 1$, then it has an infinite sequence of conservation laws of
distinct weights.}

\vspace{1pc}
In section 1, we solve the problem of equivalence, up to point transformations, for the differential systems that locally correspond to (1).  A canonical two dimensional G bundle is attached and we find a complete set of invariants whose functional relations determine the differential system up to the admissible transformations. The structure equations on the infinite prolongation space developed in section 2 enable us to find a certain subspace, rather than a quotient space, in which every conservation law has a unique representative. We compute a rough normal form of a conservation law, and this leads to the universal integrability condition to admit any higher order conservation laws. The representation of the structure group G on the conservation laws gives the aforementioned notion of \emph{weight}. In the final section, we compute the conservation laws of two classes of differential systems using the general theory. In particular, the \emph{$k_1$ flow}, the flow of the curve in the plane by the derivative of its curvature with respect to the arc length, turns out to have the \emph{KdV} property, i.e., an infinite sequence of conservation laws. 

We would like to thank Prof. Robert Bryant for his guidance and support throughout this work.

\section{Structure equations on F $\to$ M} 

Consider the following 1 forms on  $( t, x, u, p, q )$ space, where $p$ and $q$ stand for $u_x$ and $u_{xx}$ respectively.
\begin{align}
w^2    &= f(x,u,p) dt   \notag  \\
\eta_0 &= f(x,u,p) ( du - p \, dx)  \notag \\
w^1    &= dx   \notag \\
\eta_1 &= f(x,u,p) ( dp - q \, dx)  \notag  \\
\gamma_2 &= f(x,u,p) dq + g(x,u,p,q) dx \notag
\end{align}
Then the integral manifolds of 
\begin{equation} 
\mathcal{I} = \{ \, \eta_0 \wedge w^1 + \gamma_2 \wedge w^2, \, \, \eta_1 \wedge w^2, \, \, \eta_0 \wedge w^2 \, \} \, \cup \, \{ \, \eta_1 \wedge \eta_0 \, \}   
\end{equation}
on which $ \, w^1 \wedge w^2 \ne 0 \, $ correspond to the solutions of the pde  (1). Note the 1 forms introduced above satisfy the following structure equations.
\begin{align} 
dw^2 &\equiv 0 \quad \mod \quad w^2 \\
dw^1 &\equiv 0 \quad \mod \quad w^1, \, \eta_0 \notag \\
d\eta_0 &\equiv - \eta_1 \wedge w^1   \quad \mod \quad \eta_0  \notag \\
d\eta_1 &\equiv - \gamma_2 \wedge w^1 \quad \mod \quad \eta_0, \, \eta_1,   \notag
\end{align}
and if we put 
$T$ = $\frac{\partial}{\partial{t}}$,
\begin{align} \label{A:con2}
        T \lrcorner w^2    &\ne 0  \\
        T \lrcorner \eta_0 &= T \lrcorner w^1 = T \lrcorner \eta_1= T \lrcorner \gamma_2= 0  \notag \\
        T \lrcorner d\eta_0 &= T \lrcorner dw^1 = T \lrcorner d\eta_1 = T \lrcorner d\gamma_2= 0.    \notag    
\end{align}
Conversely, if a coframe $ \{ \, w^2, \eta_0, w^1, \eta_1, \gamma_2 \, \} $ of a five manifold M satisfies the conditions (3) and (4) for a nonzero vector field $T$, then the differential ideal
(2) corresponds (locally) to the pde 
\begin{equation} 
  u_t = f(x, u, p, q) u_{xxx} + g(x, u, p, q) , \notag 
\end{equation}
which includes (1) as a special case. We call such differential systems  
\emph{  quasilinear } \emph{time-independent} \emph{ third order} \emph{evolutionary} \emph{ differential systems}.  

Given such a differential system $\mathcal{I}$ on M, M has a $G_0$ $\subset$ Gl(5,R) structure, 
where $G_0$ is a subgroup whose induced action on $\bigwedge^{2}T^{*}$M preserves the 
subspace defined by $\mathcal{I}$. Equivalently, the principal right Gl(5,R) bundle over M can be 
reduced to a $G_0$ bundle via $\mathcal{I}$. Let $F_0$ denote the reduced 
$G_0$ bundle. A direct computation shows the elements of $G_0$ act on $\{$ $w^2,$ $\eta_0,$ $w^1,$ $\eta_1,$ $\gamma_2$ $\}$ by
\[
   \begin{pmatrix}  
       b^3 & .  & . & . & .\\
        . & a & . & . & .\\
        . & e & b & . & .\\
        . & h & . & \frac{a}{b} & .\\
        . & l & . & m & \frac{a}{b^2}
   \end{pmatrix}
   \begin{pmatrix}  
       w^2\\
       \eta_0\\
       w^1\\
       \eta_1\\
       \gamma_2
   \end{pmatrix},
\]
where '.' denotes 0. Note that if we allow the equivalence up to full contact transformations, $ \eta_1 $ term can be
added to $ w^1 $. 
Thus on $F_0$, we have the following structure equations.
\[  
   d \begin{pmatrix}
       \omega^2\\ \theta_0\\ \omega^1\\ \theta_1\\  s_2
     \end{pmatrix}
   = - \begin{pmatrix}
         3 \beta & . & . & . & . \\        
             .   & \alpha & . & . & . \\
             .   & \phi & \beta & . & . \\
             .   & \phi_1 & . & \alpha - \beta & . \\
             .   & \psi & . & \psi_1 & \alpha - 2 \beta
       \end{pmatrix}
   \wedge \begin{pmatrix}
             \omega^2\\ \theta_0\\ \omega^1\\ \theta_1\\  s_2
          \end{pmatrix} \]
 \[  \qquad + \begin{pmatrix}
           0 \\
        - \theta_1 \wedge \omega^1 + \theta_0 \wedge \epsilon_0\\
           \Gamma^1 \\ 
        - s_2 \wedge \omega^1 + \theta_0 \wedge \epsilon_1
                              + \theta_1 \wedge \epsilon_2   \\

           \Gamma_2
     \end{pmatrix},
\]  
with
  \begin{align}
     d\beta &= 0 \notag \\
     d\alpha &\equiv - \phi_1 \wedge \omega^1 + \phi \wedge \theta_1 \mod \theta_0. \notag
  \end{align}
Here  \: $\omega^2, \, \theta_0, \, \omega^1, \, \theta_1,  \, s_2$  \: are tautological
1 forms, \: $\alpha, \, \beta, \, \phi, \, \phi_1, \, \psi, \, \psi_1 $ \: are 
pseudo connection forms, and
  \begin{align} 
      &\epsilon_0, \, \epsilon_1, \, \epsilon_2 \, \equiv 0 \mod \, \theta_0, 
                                                     \omega^1, \theta_1, s_2 \notag\\
      &\Gamma^1, \, \Gamma_2 \qquad \, \qquad  \mbox{quadratic in } \,  \theta_0, 
                                                     \omega^1, \theta_1, s_2  \notag 
  \end{align}
representing the torsion of the pseudo connection defined by \: $\alpha, \, \beta, \, \phi,$ $\, \phi_1,$ $\, \psi,$ $\, \psi_1$.
This pseudo connection is not uniquely defined. Modifying the pseudo connection forms by
the tautological 1 forms $ \, \theta_0, \, \omega^1, \, \theta_1,  \, s_2$ , 
the torsion can be arranged as follows.
\begin{align}
   \epsilon_0    &= \epsilon_1 = 0 \notag \\
   \Gamma^1    &= a_1 \theta_1 \wedge \omega^1 + R s_2 \wedge \omega^1 \notag \\ 
   \epsilon_2 &= a_2 \omega^1 + R s_2 \notag \\
   \Gamma_2    &= a_3 s_2 \wedge \omega^1 \notag ,
\end{align}
where $a_1, a_2, a_3, R$ are functions on $F_0$. Note that, again, since the equivalence up to point transformations is allowed, the differential ideal generated by $\{ \, \theta_0, \, \omega^1 \, \}$ is integrable. Also, the original differential system corresponds (locally) to the
equations of the form
\[
   u_t \, = \, f(\, x, u, p \, ) u_{xxx} + g(\, x, u, p, q \, )
\]
iff $ R = 0 $, which we assume from now on. 

At this stage, no more absorption is possible without changing the form of the torsion given above, i.e., the remaining torsion terms are 
unabsorbable. By taking the exterior derivative of the structure equations with the reduced torsion, we get
\begin{align}
     da_1  &\equiv \phi + a_1 ( \alpha - \beta ) 
                               \mod  \theta_0, \, \omega^1, \, \theta_1   \notag \\
     da_2  &\equiv 2 \phi_1 - \psi_1 + a_2 \beta  
                               \mod  \theta_0, \, \omega^1, \, \theta_1, \,  s_2\notag \\
     da_3  &\equiv \phi_1 + \psi_1 +   a_3 \beta  \notag 
                                  \mod  \theta_0, \, \omega^1, \, \theta_1, \,  s_2 .
\end{align}
Hence, the $G_0$ structure can be reduced to a $G_1 \subset G_0$ structure by requiring
$\, a_1 = 0, \, a_2 = 0, \, a_3 = 0$. In other words, we restrict to the subbundle $F_1$ 
of $F_0$ cut out by these relations. 
The structure equations on $F_1$ now become
\begin{equation}
   d \begin{pmatrix}
       \omega^2\\ \theta_0\\ \omega^1\\ \theta_1\\  s_2
     \end{pmatrix}
 \notag
   = - \begin{pmatrix}
         3 \beta & . & . & . & . \\        
             .   & \alpha & . & . & . \\
             .   &  . & \beta & . & . \\
             .   &  . & . & \alpha - \beta & . \\
             .   & \psi & . & . & \alpha - 2 \beta
       \end{pmatrix}
 \notag
   \wedge \notag \begin{pmatrix}
             \omega^2\\ \theta_0\\ \omega^1\\ \theta_1\\  s_2
          \end{pmatrix}
 \notag
\end{equation}
\begin{equation}
   + \begin{pmatrix}
           0 \\
        - \theta_1 \wedge \omega^1 \\
          \theta_0 \wedge \phi \\ 
        - s_2 \wedge \omega^1 + \theta_0 \wedge \phi_1 \\
          \theta_1 \wedge \psi_1
     \end{pmatrix}, \notag
\end{equation}
where $ \phi \equiv 0 \mod  \theta_0, \omega^1, \theta_1 \,$ and
  $\phi_1, \psi_1 \equiv 0 \mod  \theta_0, \omega^1, \theta_1,  s_2. \,$
Again, by modifying $\alpha,  \psi$, we can arrange so that
\begin{align}
      \phi &= E \omega^1 + F \theta_1 \notag \\
      \phi_1 &= a_4 \omega^1 + K s_2 \notag \\ 
      \psi_1 &= G \omega^1 + 2 K s_2 \notag ,
\end{align}
where $ \, E, F, G, K, a_4 \, $ are now functions on $F_1$ that represent the unabsorbable torsion. Taking the exterior derivative again of the structure equations with this modified torsion gives, among
other things, that $ a_4 $ can be translated to 0 by the group action corresponding to $ \psi $. Thus 
we can reduce the $G_1$ structure $F_1$ to a two dimensional 
$\mbox{G} \subset G_1$ structure $\mbox{F}$ on which $ a_4 = 0$, and $ \psi \equiv 0 \mod
 \theta_0, \, \omega^1, \, \theta_1, \,  s_2.$  The elements of $\mbox{G}$ are of the form
\[
   \begin{pmatrix}  
       b^3 & .  & . & . & .\\
        . & a & . & . & .\\
        . & . & b & . & .\\
        . & . & . & \frac{a}{b} & .\\
        . & . & . & . & \frac{a}{b^2}
   \end{pmatrix}.
\]

Summarizing the results, we have the following proposition.

\begin{proposition} \label{P:streq}
 Let \textnormal{M} be a five manifold with quasilinear time-independent third order evolutionary
 differential system $\mathcal{I}$ that locally corresponds to the equation of the form
\begin{center}
    $u_t \, = \, f(\, x, u, p \, ) u_{xxx} + g(\, x, u, p, q \, ).$
\end{center}
 Then there is an induced \textnormal{G} $\subset$ \textnormal{Gl(5,R)} 
bundle \textnormal{F} $\to$ \textnormal{M} with the structure equations
\[
    d \begin{pmatrix}
           \omega^2\\ \theta_0\\ \omega^1\\ \theta_1\\  s_2
      \end{pmatrix}
     = - \begin{pmatrix}
           3 \beta & . & . & . & . \\        
               .   & \alpha & . & . & . \\
               .   & . & \beta & . & . \\
               .   & . & . & \alpha - \beta & . \\
               .   & . & . & . & \alpha - 2 \beta
         \end{pmatrix}
   \wedge \begin{pmatrix}
             \omega^2\\ \theta_0\\ \omega^1\\ \theta_1\\  s_2
          \end{pmatrix}
\]
\[
   \; \: \quad  \quad \qquad \qquad  \quad  + \begin{pmatrix}
             0 \\
          - \theta_1 \wedge \omega^1 \\
            \theta_0 \wedge ( E \omega^1 +  F \theta_1  ) \\ 
          - s_2 \wedge \omega^1 + K \theta_0 \wedge s_2\\
            \theta_1 \wedge ( G \omega^1 + 2 K s_2 ) + \theta_0 \wedge ( L \omega^1 +  
                                                       M \theta_1 + N s_2 )            
         \end{pmatrix}
\]
    \begin{align}
           d\beta  &= 0 \notag \\
           d\alpha  &= - E \, \theta_1 \wedge \omega^1 - K \, s_2 \wedge \omega^1 + \theta_0 \wedge \epsilon, \notag
   \end{align}
where $\epsilon  \equiv 0 \mod \omega^1, \, \theta_1, \, s_2.$ The complete set of invariants of this \textnormal{G} structure consists of E, F, G, K, L, M, N, $\epsilon $ and their successive covariant derivatives. If all of the invariants are 0, the differential system is equivalent to the system generated by $u_t = u_{xxx}$ via a point transformation, i.e., a diffeomorphism of \textnormal{M} that preserves the ideal generated by $\, \{ \, \theta_0, \, \omega^1 \, \}$. By construction, if $\mathcal{I}$ and $\mathcal{I}'$ are two differential systems such that $\varphi^* \mathcal{ I }' = \mathcal{ I }$ for a point transformation $\varphi$ of \textnormal{M}, then the induced \textnormal{G} structures \textnormal{F} and \textnormal{F}$'$ are equivalent. \\
\end{proposition} 

Given an evolutionary pde  (1),
we mention that 
  $ K = 0 $ 
iff the pde is of the form
\begin{center} 
    $u_t = ( f(x, u, p) q )_x + g(x, u, p) q + h(x, u, p).$
\end{center}
In this case, $ E = -N,$ and 
$ N = 0 $ iff
  \begin{align}
     4 f_{p}^2 &= 3 f f_{pp}, \notag \\
     4 ( \, f_p f_x + p f_p f_u \, ) &= 3 f ( \, f_{px} + p f_{pu} - f_u ). \notag
  \end{align}

\section{Structure equations on F$_{\infty} \rightarrow$ F $\rightarrow$ M} 

Consider the linear pde
\[ u_t = u_{xxx}. \]
Since every $t$ derivative of $u$ can be replaced by $xxx$ derivative for any solution $u$,
the k + 3 th jet space $J_{k + 3}$ of $u$ has as a coordinate system
\[ \{ \, x, \, t, \, u = p_0, \, p_1, \, p_2, \, p_3, \,  . . . . p_{k + 3} \,   \}  , \]
where $p_1 = u_x, \, p_2 = u_{xx}$, and $p_{i} = u_{xxx . . . x}$. In terms of differential forms,  we introduce the following 1 forms on $J_{k + 3}$.
\begin{align}
               w^1 &= dx, \, w^2 = dt,                 \notag  \\
       \gamma_i       &=dp_i   \; \qquad \qquad \qquad \qquad \quad  \mbox{for} \; \; 0 \leq i \leq k + 3   , \notag    \\
      \eta_i          &= \gamma_i - p_{i + 1} dx \qquad \qquad \quad \, \, \, \,  \mbox{for} \; \; 0 \leq i \leq k + 2 , \notag  \\
       \phi_i          &= \eta_i \qquad \qquad \quad - p_{i + 3} dt  \quad   \mbox{for} \; \; 0 \leq i \leq k  .\notag    
\end{align}
Then $ \, w^1, \, w^2, \, \phi_0, \, \phi_1, \, . . . \phi_k, \, \eta_{k + 1}, \, \eta_{k + 2}, \,  \gamma_{k +3} \,$ form a basis of $T^*J_{k + 3}$ and, in fact,  they satisfy the following structure equations.
\begin{align}
   - d\eta_i       &= \eta_{i + 1} \wedge w^1 = \gamma_{i + 1} \wedge w^1 \quad \mbox{for} \; \; 0   \leq i \leq k + 1  \notag \\
   - d\eta_{k + 2} &= \gamma_{k + 3} \wedge w^1   \notag  
\end{align}
\begin{align}                   
   - d\phi_i       &= \phi_{i + 1} \wedge w^1 + \phi_{i + 3} \wedge w^2 \notag\\
                   &= \eta_{i + 1} \wedge w^1 + \gamma_{i + 3} \wedge w^2  \quad \mbox{for} \; \; 0   \leq i \leq k - 3  \notag \\
   - d\phi_{k - 2} &=  \phi_{k - 1} \wedge w^1 + \eta_{k + 1} \wedge w^2    \notag   \\ 
   - d\phi_{k - 1} &=  \, \, \phi_k \, \wedge \, w^1 + \eta_{k + 2} \wedge w^2    \notag   \\ 
   - d\phi_k       &=  \eta_{k + 1} \wedge w^1 +  \gamma_{k + 3} \wedge w^2 .   \notag 
\end{align}
Note the original differential system for  $u_t = u_{xxx}$ on $(\, x, t, u, p_1, p_2 \,)$ space is generated by
\begin{align}
    \: &\{ \: \eta_0  \wedge w^1 + \gamma_2 \wedge w^2 \: \}         \notag \\
                  \cup \, &\{ \: \eta_1  \wedge w^2,  \:   \eta_0 \wedge w^2 \: \}  \notag  \\
                  \cup \, &\{ \: \eta_1 \wedge \eta_0 \: \}. \notag
\end{align}   
Moreover, the structure equations show, with respect to the given basis, the k + 3 th 
prolongation of this differential system on $J_{k + 3}$ is generated by
\begin{align}
    \: &\{ \: \phi_0, \, \phi_1, \, \phi_2, \,  . . . \phi_k \: \} \notag \\
            \cup \,  &\{ \: \eta_{k + 1} \wedge w^1 + \gamma_{k + 3} \wedge w^2 \: \}   \notag \\
            \cup  \,  &\{ \: \eta_{k + 2} \wedge w^2 , \:  \eta_{k + 1} \wedge w^2 \: \} \notag \\
            \cup \, &\{ \: \eta_{k + 2} \wedge \eta_{k + 1} \: \}. \notag
\end{align}

For a general differential system $\mathcal{I}$ on M, let M$_{k + 3}$ denote the k + 3 th prolongation space, M$_{\infty}$ = $\lim_{k \to \infty}$ M$_{k + 3}$, and let $p_3, p_4, p_5$, ... be the prolongation variables. Also define F$_{k + 3}$ \, and  \, F$_{\infty}$ as the pullback of  
$\,$ M$_{k + 3}$ and $\,$  M$_{\infty}$ via $\pi$ respectively, where $\pi$ is the projection map 
$\pi :$ F $\to$ M of the G structure defined earlier. Motivated by the linear example above, we observe the following.

\begin{theorem} 
There is a sequence of time-independent functions $r_{k + 3}$ on \textnormal{F}$_{k + 2}$  $\:$ for $\, \, k \geq 2$ such that if we define 
\begin{align}
   s_i &= \, dp_i + p_i ( \, \alpha - i\beta \, )  \quad \textnormal{for} \, \,  \, i \geq 3, \notag \\
  \sigma_2  &= s_2, \quad \sigma_3  = s_3, \quad  \sigma_4  = s_4 + p_3 K s_2 \notag \\ 
        &\notag  \\
  \sigma_{k + 3} &= s_{k + 3} + r_{k + 3} \omega^1  \quad \; \; \; \; \;   \textnormal{for} \, \,  k \geq 2,  \notag      \\
  \theta_{k} &= \sigma_{k} - p_{k + 1} w^1  \qquad \qquad \quad \textnormal{for} \, \, \,   k \geq 2 , \notag   \\
  \pi_{k} &= \theta_{k} \quad \quad \qquad - p_{k + 3} \omega^2  \quad \: \textnormal{for} \, \, \,  k \geq 0 , \notag  
\end{align}
then 

 $1.$ \:  $ \{ \,  \alpha, \, \beta, \, \omega^1, \, \omega^2, \, \pi_0, \, \pi_1, \, . . . \, 
           \pi_k, \, \theta_{k + 1}, \, \theta_{k + 2}, \, \sigma_{k +3} \, \} 
           \,$ is a basis of $  T^*$\textnormal{F}$_{k + 3}$   for $\, k \geq 0.  $ 

 $2.$ For  $\, k \geq 0,  $ \begin{align}
   - d\pi_i  &\equiv \pi_{i + 1} \wedge \omega^1 + \pi_{i + 3} \wedge \omega^2     
              \: \mod \: \pi_0 \, , .  .  . \, \pi_{i} \, \qquad \, \: \textnormal{for} \, \, \, 0 \: \leq \:  i \: \leq k-3 \notag \\
   - d\pi_{k - 2}    &\equiv \pi_{k - 1} \wedge \omega^1 + \theta_{k + 1} \wedge \omega^2 
               \mod \: \pi_0  \, ,  .  .  . \,  \pi_{k - 2}      \notag   \\ 
   - d\pi_{k - 1} &\equiv  \pi_k \, \wedge \, \omega^1 + \theta_{k + 2} \wedge \omega^2  
              \; \: \: \mod \: \pi_0  \, ,  .  .  .  \,  \pi_{k - 1} \,   \notag   \\
   - d\pi_k   &\equiv \theta_{k + 1} \wedge \omega^1 +  \sigma_{k + 3} \wedge \omega^2  
               \: \: \mod \: \pi_0  \, ,  .  .  .   \, \pi_{k}   . \notag 
  \end{align}

 $3$. For $\, i \geq 0,$ 
  \begin{align}
    -d\theta_i &\equiv  ( \, \alpha - i \beta \, ) \wedge \theta_i + \theta_{i + 1} \wedge  
   \omega^1 \quad \mod \quad \theta_0, \, \theta_1, \, . . . \,  \theta_{i - 2}, \, \theta_{i - 1}.  \notag 
  \end{align} 

Note that $\, 2.$ implies the k + 3 th prolonged ideal $\mathcal{I}_{k + 3}$ on \textnormal{F}$_{k + 3}$ is generated by
\begin{align}  \mathcal{I}_{k + 3} \: = \: 
   \quad &\{ \: \pi_0, \, \pi_1, \, \pi_2, \,  . . . \pi_k \: \} \notag \\
   \cup \,  &\{ \: \theta_{k + 1} \wedge \omega^1 + \sigma_{k + 3} \wedge \omega^2 \: \}   \notag \\
  \cup \,  &\{ \: \theta_{k + 2} \wedge \omega^2 , \:  \theta_{k + 1} \wedge \omega^2 \: \} \notag  \\
  \cup \,  &\{ \: \theta_{k + 2} \wedge \theta_{k + 1} \: \}. \notag     
\end{align}    
Also by definition, 
  \begin{align}  
    \pi_{i} \wedge \omega^1 + \theta_{i + 2} \wedge \omega^2 &=
    \theta_{i} \wedge \omega^1 + \sigma_{i + 2} \wedge \omega^2 \notag \\
    \theta_{i} \wedge \omega^1 &= \sigma_{i} \wedge \omega^1 \notag \\
    \pi_{i} \wedge \omega^2  &=  \theta_{i} \wedge \omega^2 . \notag
  \end{align}
\end{theorem}

\noindent \emph{Proof}

1. It follows from the definition, as long as $r_{k + 3} \in$ $C^{\infty}$(F$_{k + 2}$), which is to be checked below. 

2. First, let's consider the original differential system 
$\mathcal{I}$ = $\mathcal{I}_2$ on F = F$_2$. ( For simplicity, we use the same letter 
$\mathcal{I}$ to denote $\pi^*$ $\mathcal{I}$ on F, where $\pi$ : F $\to$ M. )
\begin{align}  \mathcal{I} = \: \mathcal{I}_{2} = \: 
 \,  &\{ \: \theta_{0} \wedge \omega^1 + \sigma_{2} \wedge \omega^2 \: \}   \notag \\
  \cup \,  &\{ \: \theta_{1} \wedge \omega^2 , \:  \theta_{0} \wedge \omega^2 \: \} \notag \\
  \cup \, &\{ \: \theta_1 \wedge \theta_0 \: \}. \notag
\end{align} 
Since $\:  \theta_{0} \wedge \omega^2 \, \in \mathcal{I}$ , $\theta_0$ is a multiple of
$\omega^2$ on every integral manifold of $\mathcal{I}$ on which $\omega^1 \wedge \omega^2 \ne$ 0.
This suggests introducing 
\begin{align}
              \textnormal{F}_3 &= \textnormal{F} \times \{\, p_3 \, \} \notag \\
             \pi_0 &= \theta_0 - p_3 \omega^2. \notag 
\end{align}
But,
\[
   \theta_{0} \wedge \omega^1 + s_{2} \wedge \omega^2 \,
       \equiv   \theta_2 \wedge \omega^2   \quad  \, \mod \; \pi_0 
\]
by definition of $\theta_2$. Also from the  structure equations on F,
\[ 
    - d\pi_0 = - d( \, \theta_0 - p_3 \omega^2 \, ) \equiv \theta_1 \wedge \omega^1 + \sigma_3 \wedge \omega^2  \quad \mod  \, \, \, \pi_0 . 
\]
This implies $\mathcal{I}_3$ is generated by
\begin{align}  \mathcal{I}_{3} = \: 
   \quad &\{ \: \pi_0 \: \} \notag \\
   \cup \,  &\{ \: \theta_{1} \wedge \omega^1 + \sigma_{3} \wedge \omega^2 \: \}   \notag \\
  \cup \,  &\{ \: \theta_{2} \wedge \omega^2 , \:  \theta_{1} \wedge \omega^2 \: \} \notag \\
  \cup \, &\{ \: \theta_2 \wedge \theta_1 \: \}, \notag
\end{align}  
which is the \emph{$p_3$ prolongation}  $\mathcal{I}_3$ of $\mathcal{I}$. From the arguments above, it is clear that every integral manifold of $\mathcal{I}_3$ in F$_3$ on which $\omega^1 \wedge \omega^2 \ne$ 0 is a graph over an integral manifold of $\mathcal{I}$ in F. Note that 
\[ \{ \,  \alpha, \, \beta, \, \omega^1, \, \omega^2, \, \pi_0, \, \theta_{1}, \, \theta_{2}, 
       \, \sigma_{3} \, \}  \]
is a basis of $T^*$F$_3$.

Next, we define the \emph{$p_4$ prolongation} by the similar procedure, introducing
\begin{align}
               \textnormal{F}_4 &= \textnormal{F}_3 \times \{\, p_4 \, \} \notag \\
             \pi_1 &= \theta_1 - p_4 \omega^2 . \notag 
     \end{align}
Again,
\[
   \theta_{1} \wedge \omega^1 + \sigma_{3} \wedge \omega^2 \,
       \equiv   \theta_3 \wedge \omega^2   \quad  \, \mod \; \pi_0, \, \pi_1
\]
and
\[ 
    - d\pi_1 = - d( \, \theta_1 - p_4 \omega^2 \, ) \equiv \,
                \theta_2 \wedge \omega^1 + \sigma_4 \wedge \omega^2  \quad \mod  \; \pi_0 , \, \pi_1 . 
\]
Hence,
\begin{align}  \mathcal{I}_{4} = \: 
   \quad &\{ \: \pi_0 , \, \pi_1 \: \} \notag \\
   \cup \,  &\{ \: \theta_{2} \wedge \omega^1 + \sigma_{4} \wedge \omega^2 \: \}   \notag \\
  \cup \,  &\{ \: \theta_{3} \wedge \omega^2 , \:  \theta_{2} \wedge \omega^2 \: \} \notag \\
   \cup \, &\{ \: \theta_3 \wedge \theta_2 \: \}. \notag 
\end{align}

The claim for \emph{$p_5, \, p_6, \, p_7$ prolongations} can also be checked as above by direct computations. 

Inductively, suppose there exists a sequence of time independent functions $r_5, \,$ $r_6, \,$ $. . . \,$ $r_{k + 3}$, $k \geq 4, \,$ 
such that the structure equations in $2.$ of the theorem is true. It suffices to show there exists time independent  
$r_{k + 4}$ $\in$ $C^{\infty}($F$_{k + 3}$) such that 
\[
    - d\pi_{k + 1}   \equiv  \theta_{k + 2} \wedge \omega^1  +  \sigma_{k + 4} \wedge \omega^2  
               \mod \: \pi_0  , \,  .  .  .  . \, \pi_{k + 1} .
\]
We only verify the claim $\mod \alpha, \, \beta \; $ for simplicity. The terms containing  $\, \alpha \, $ and  $\, \beta \, $ can be checked by counting the weights of the representation of the structural group G on F$_{\infty}$. For brevity, we use 
 $\{ \,  p_i \, \}$ to denote a time independent function on F$_i$, i.e., it's 
a function of $p_3$,  $p_4$,  $...$   $p_i$, and the complete set of invariants 
of the original G structure F $\to$ M.  

Notice from the structure equation on F $ \to $ M,
\begin{align}
    d\omega^2 &\equiv  0   \notag \\
    d\theta_0 &\equiv  p_4 \, \omega^1 \wedge \omega^2 \notag \\
    d\omega^1 &\equiv  \{ \, p_3  \, \} \omega^1 \wedge \omega^2 \notag \\
    d\theta_1 &\equiv  \{ \, p_5  \, \}  \omega^1 \wedge \omega^2 \notag \\
    d s_2     &\equiv  \{ \, p_4 \, \}  \omega^1 \wedge \omega^2 \notag \\
    d\alpha   &\equiv  \{ \, p_5  \, \}  \omega^1 \wedge \omega^2 
                    \mod \alpha, \, \beta, \, \pi_0, \, \pi_1, \, \pi_2.   \notag 
\end{align}
Thus,
\begin{align}
   d\pi_{k + 1} &= d( \, \theta_{k +1} - p_{k + 4} \omega^2 \, ) \notag \\
                &= d( \, dp_{k + 1} +p_{k + 1} ( \, \alpha - (k + 1) \beta \, ) 
                        + r_{k + 1} \omega^1 - p_{k + 2} \omega^1 - p_{k + 4} \omega^2 \, ) \notag \\
        &\equiv p_{k + 1} \{ \, p_5 \, \}  \omega^1 \wedge \omega^2  
     + dr_{k + 1} \wedge  \omega^1 + r_{k + 1} \{ \, p_3  \, \}  \omega^1 \wedge \omega^2 
                  \notag \\
      &\quad + p_{k + 2} \{ \, p_3 \, \}  \omega^1 \wedge \omega^2
                  - dp_{k + 2} \wedge \omega^1 - dp_{k + 4} \wedge \omega^2 \notag \\
     &\equiv  - dp_{k + 2} \wedge \omega^1 - dp_{k + 4} \wedge \omega^2 +  dr_{k + 1} \wedge  \omega^1 
                         + \{  \, p_{k + 3} \, \} \,  \omega^1 \wedge \omega^2 \notag \\
    &  \mod \alpha, \, \beta, \,  \pi_0, \, \pi_1, \, \pi_2 \notag
\end{align}
 since  $r_{k + 1} \in$ $C^{\infty}($F$_{k}$) is time independent by induction hypothesis and $k \geq 4$.
Also,
\[
  \theta_{k + 2} \wedge \omega^1 + \sigma_{k + 4} \wedge \omega^2 \equiv \,
           dp_{k + 2} \wedge \omega^1 + dp_{k + 4} \wedge \omega^2 + r_{k + 4} \omega^1 \wedge \omega^2
                 \mod \alpha, \, \beta ,
\]
by definition of $r_{k + 4}$. Now, it is sufficient to show 
\[ dr_{k + 1} \wedge \omega^1 \equiv \{ \, p_{k + 3} \, \} \,  \omega^1 \wedge \omega^2  \mod  \, \alpha, \, \beta, \, \pi_0, \, . . . \, \pi_k.
\]
But, again, $r_{k + 1} \in$ $C^{\infty}($F$_{k}$) is time independent and by definition of the 1 form $\pi_i$,
\[ dp_i \equiv \, p_{i + 3} \omega^2 \, \mod  \, \alpha, \, \beta, \, \omega^1,  \, \pi_i.
\]
In particular,
\[ dp_k \equiv \, p_{k + 3} \omega^2 \, \mod  \, \alpha, \, \beta, \, \omega^1,  \, \pi_k,
\]
and $r_{k + 4} \in$ $C^{\infty}($F$_{k+3}$) is uniquely determined. 

3. Straightforward computations. We omit the proof. 
\begin{flushright} $\square$ \end{flushright}

 \section{Conservation laws}

In this section, we compute a rough normal form of conservation laws together with
a universal integrability condition to admit any higher order conservation laws
using the structure equations on F$_{\infty} \to$ F. 

\subsection{$\bigwedge^2$ $\mathcal{I}_{\infty}$ $/ d\mathcal{I}_{\infty}$}

Following [BG1], we start with a definition.
\begin{definition} 
A conservation law is a class \textnormal{[} $\Phi$ \textnormal{]} $\in$
$\bigwedge^{2}$ $\mathcal{I}_{\infty}$ / $d\mathcal{I}_{\infty}$ generated by a closed 
2 form $\Phi$. The space of conservation laws will be denoted by $\mathcal{C}_{\infty}$.
\end{definition}
\noindent Here $d\mathcal{I}_{\infty}$ is the subspace of exact 2 forms in the ideal
 of the form $d(\: \lambda_0 \pi_0 \, + \, \lambda_1 \pi_1 \, + \,  . . . \, 
\lambda_m \pi_m \: )$ for some functions  $\lambda_i'$s on F$_{\infty}$. 
For such a closed 2 form $\Phi$, consider a 1 form $\theta$ such that
\[ d\theta \, = \, \Phi  .\]
If N is an integral manifold of $\mathcal{I}$, hence of $\mathcal{I}_{\infty}$ , it follows 
by Stokes theorem,
\[ 0 \, = \, \int_{N} \Phi \, = \,  \int_{N} d\theta \, = \, \int_{\partial N} \theta ,\]
which becomes the classical conservation law for an evolutionary pde once we impose
some nice decay conditions on the solutions. [BG1] shows that, at least locally,
the cohomology  $\mathcal{C_{\infty}}$ defined as above is isomorphic to the classical space of
conservation laws
\begin{center}
    $\mathcal{C}_{\infty} \, \cong \{ \: \theta \in \bigwedge^1 \textnormal{F}_{\infty} \: \mid  \:
                              d\theta \equiv 0 \, \mod \mathcal{I}_{\infty} \: \} \: / \: \{ \, \bigwedge^1 \mathcal{I}_{\infty} \, \}  \cup   \{ \: df \: \mid \: f \in C^{\infty}( \textnormal{F}_{\infty} ) \: \}$. \\
\end{center} 

Instead of dealing directly with the quotient $\bigwedge^{2}$ $\mathcal{I}_{\infty}$ / $d\mathcal{I}_{\infty}$, it is convenient to introduce a section of the projection
\begin{center}
  $\bigwedge^{2}  \mathcal{I}_{\infty} \to \bigwedge^{2} \mathcal{I}_{\infty} / d\mathcal{I}_{\infty}.$
\end{center}

\begin{theorem}
There exists a sequence of subspaces $\, \mathcal{H}_{k + 3} \subset 
      \bigwedge^{2} \mathcal{I}_{k + 3} \;$
      for $\, k \geq  \, 0 \,$  such that 
\begin{align}
 1. \quad \; \mathcal{H}_{k + 3} &= \, \, \mathcal{H}_{k + 2} \, \cup \{ \, \theta_{k + 2} \wedge \theta_{k + 1} \, \}
       \, \cup \, \{ \, \theta_{k + 2} \wedge \theta_{k} \, \}  \notag \\
         & \qquad \qquad \, \cup \, \{ \, \sigma_{k + 3} \wedge \pi_i \mid \, 0 \, \leq \, i \, \leq \, k \, \} \notag \\
         & \qquad \qquad  \, \cup \, \{ \, \theta_{k + 2} \wedge \omega^2 \, \}.  \notag  \\
     \mathcal{H}_{k + 3} \, &\cong \, \textnormal{$\bigwedge$}^2 \mathcal{I}_{k + 3} /   
                 d\mathcal{I}_{k + 3}  \quad  \textnormal{for} \: \: k \, \geq \, 0  .\notag  \\
 2. \quad   Set  \quad &\; \mathcal{H}_{\infty} \, = \, \lim_{k \to \infty} 
           \mathcal{H}_{k + 3} \, \quad then, \notag \\
    \mathcal{H}_{\infty} \, &\cong \, \textnormal{$\bigwedge$}^2  \mathcal{I}_{\infty}  / d\mathcal{I}_{\infty}. \notag
\end{align}
\noindent We  define
\begin{align}
  \mathcal{H}_1 \, = \, & \{ \, \theta_0 \wedge \omega^1 + \sigma_2 \wedge \omega^2 \, \} \notag \\
               \cup  \, & \{ \, \theta_1 \wedge \omega^2 , \, \theta_0 \wedge \omega^2 \, \} , \notag \\
  \mathcal{H}_2 \, = \, & \mathcal{H}_1 \, \cup \, \{ \, \theta_1 \wedge \theta_0 \, \}. \notag
\end{align}
\end{theorem}

\begin{corollary}
 Every conservation law has a unique representative in $\mathcal{H}_{\infty}$. \end{corollary}

\vspace{1pc}
\noindent \emph{Proof of the Theorem}

1. Suppose $\, \Phi \in \bigwedge^2 \mathcal{I}_{k + 3} \, $  for some k $\geq$ 0. Then
\begin{align}
  \Phi &= \pi_0 \wedge \epsilon_0 + \pi_1 \wedge \epsilon_1 + \,  . . . \, \pi_k \wedge \epsilon_k 
        \, + \sum_{0}^{k} c_{ij} \, \pi_i \wedge \pi_j   \notag \\
       &\, + v \, ( \: \theta_{k + 1} \wedge \omega^1 + \sigma_{k + 3} \wedge \omega^2 \: )   \notag \\
       &\, + y \, \theta_{k + 2} \wedge \omega^2  \: + z \, \theta_{k + 1} \wedge \omega^2 \:  \notag \\
       &\, + w \, \theta_{k + 2} \wedge \theta_{k + 1} \: , \notag   
\end{align}
where $\epsilon_i \equiv 0 \mod \, \omega^1, \, \omega^2, \, \theta_{k + 1} , \, \theta_{k + 2}, \, \sigma_{k + 3} \, $ and $ c_{ij} = - c_{ji}, \, v, \, y, \, z, \, w \, $ are functions on F$_{k + 3}$. But, since 
\[  - d\pi_k   \equiv   \theta_{k + 1} \wedge \omega^1 +  \sigma_{k + 3} \wedge \omega^2  
               \mod \: \pi_0,  \, \pi_1, \, .  .  .  . \, \pi_{k} , \,  
\]
 mod $d\mathcal{I}_{k + 3}$, we may take $ v = 0 $ by modifying $\epsilon_i '$s and $c_{ij} '$s if necessary. Also, since 
\[ - \pi_i \wedge \omega^1 \, \equiv \, \theta_{i + 2} \wedge \omega^2  \: 
           \mod   d\pi_{i - 1}, \,  \pi_{i - 1}, \, . . . \, \pi_0  \, \qquad \,   \mbox{for}
           \: 1 \, \leq \, i \, \leq k ,  \]
we may take $\, \epsilon_i \equiv 0 \, \mod \omega^2, \, \theta_{k + 1} , \, \theta_{k + 2}, \, \sigma_{k + 3}$ for $ 1 \leq i \leq k $ again by successively modifying $\epsilon_i'$ and $c_{ij}'$s if necessary. At this stage, 
no more modification  $ \mod d\mathcal{I}_{k + 3} $ is possible. Let's denote this subspace by 
${H}_{k + 3} \subset \bigwedge^2 \mathcal{I}_{k + 3}$ for each $k$ $\geq$ $0$, which satisfies
\begin{center}
    ${H}_{k + 3} \oplus d\mathcal{I}_{k + 3} = \bigwedge^2 \mathcal{I}_{k + 3}$. 
\end{center}

Next, we modify ${H}_{k + 3}$ by $d\mathcal{I}_{k + 3}$ to construct another subspace 
$\mathcal{H}_{k + 3} \subset \bigwedge^2 \mathcal{I}_{k + 3}$ complementary to $d\mathcal{I}_{k + 3}$,
but which satisfies
\[ \mathcal{H}_{k + 3} \subset \mathcal{H}_{k + 4}.  \]
First, note that as we prolong from $\mathcal{I}_{k + 3}$ to $\mathcal{I}_{k + 4}$, 
\begin{align} 
     \pi_i  \, &\to \, \pi_i  \quad  \mbox{for} \, \, 0 \leq i \leq k  \notag \\
     \theta_{k + 1} \,  &\to \,  \pi_{k + 1}  \notag \\
     \theta_{k + 2} \,  &\to \,  \theta_{k + 2}  \notag \\
     \sigma_{k + 3} \,  &\to \,  \theta_{k + 3}  \notag \\
                        &\quad \, \, \, \,  \sigma_{k + 4}.   \notag
\end{align}
This observation together with the fact
\begin{align}
    \pi_i \wedge \pi_{k + 1} &\equiv \pi_i \wedge \theta_{k + 1}  \mod  \pi_i \wedge \omega^2  \notag \\
    \pi_i \wedge \theta_{k + 3} &\equiv \pi_i \wedge \sigma_{k + 3}  \mod  \pi_i \wedge \omega^1  \notag \\
  - \pi_i \wedge \omega^1  &\equiv \theta_{i + 2} \wedge \omega^2  
     \mod d\pi_{i - 1}, \,  \pi_{i - 1}, \, . . . \, \pi_0  \, \,  \, \, \mbox{for} \, \, \,  1 \leq i \leq k  \notag  \\
    \theta_{k + 3} \wedge \pi_{k + 1} &\equiv \theta_{k + 3} \wedge \theta_{k + 1} 
                       \mod \theta_{k + 3} \wedge \omega^2 \notag 
\end{align}
implies
\begin{align}
      H_{k + 4} \equiv H_{k + 3} \, &\cup \{ \,  \theta_{k + 3} \wedge \theta_{k + 2}  \, \}
       \, \cup \, \{ \, \theta_{k + 3} \wedge \theta_{k + 1} \, \}  \notag \\
         &\cup \, \{ \, \sigma_{k + 4} \wedge \pi_i \mid \, 0 \, \leq \, i \, \leq \, k+1 \, \} \notag \\
         &\cup \, \{ \, \theta_{k + 3} \wedge \omega^2 \, \}  \notag  
             \qquad   \mod d\mathcal{I}_{k + 4}.      
\end{align}
\noindent Now, we inductively define $\mathcal{H}_{k + 4}$ by
\begin{align}
 \mathcal{H}_{k + 4} = \, \, \mathcal{H}_{k + 3} \, &\cup \{ \, \theta_{k + 3} \wedge \theta_{k + 2} \, \}
       \, \cup \, \{ \, \theta_{k + 3} \wedge \theta_{k + 1} \, \}  \notag \\
         &\cup \, \{ \, \sigma_{k + 4} \wedge \pi_i \mid \, 0 \, \leq \, i \, \leq \, k+1 \, \} \notag \\
         &\cup \, \{ \, \theta_{k + 3} \wedge \omega^2 \, \}.  \notag 
\end{align}

2. Suppose
 \begin{equation} 
d( \, \lambda_0 \pi_0 + \lambda_1 \pi_1 + \, . . . \, \lambda_m \pi_m \,) \in \mathcal{H}_{k + 3}. \notag
 \end{equation} 
There exists $l \geq k$ such that
\begin{center}
 $  \, \lambda_0 \pi_0 + \lambda_1 \pi_1 + \, . . . \, \lambda_m \pi_m \, \in \mathcal{I}_{l + 3}. $
\end{center}
But $\mathcal{H}_{k + 3}  \, \subset  \, \mathcal{H}_{l + 3}, \, \, \mathcal{H}_{l + 3} 
             \, \cap  \,     d\mathcal{I}_{l + 3}  = 0$  
implies
\[ d( \, \lambda_0 \pi_0 + \lambda_1 \pi_1 + \, . . . \, \lambda_m \pi_m \,) = 0  .\, \]
\begin{flushright} $\square$ \end{flushright}

\subsection{ Sequence of linear differential systems } 
 
Consider a conservation law $\Phi$ that belongs to  $\mathcal{H}_{k + 3}$ pointwise.
\begin{align}
 \Phi = \, \, &A ( \, \theta_0 \wedge \omega^1 + \sigma_2 \wedge \omega^2 \, ) 
                - \, B \, \theta_1 \wedge \omega^2 + \, C \, \theta_0 \wedge \omega^2 \notag \\
          + & \sum_{i = 0}^{k + 1} T_{i + 1, i} \, \theta_{i + 1} \wedge \theta_i \notag \\
          + & \sum_{i = 0}^{k} X_{i + 2, i} \, \theta_{i + 2} \wedge \theta_i 
          + \sum_{j = 0}^{k} \sum_{i = 0}^{j} X_{j + 3, i} \, \sigma_{j + 3} \wedge \pi_i \notag \\
          + & \sum_{i = 2}^{k + 2} Y_{i} \, \theta_{i} \wedge \omega^2 ,\notag
\end{align}
where $ A, B, C, T_{i + 1, i}, X_{i + 2, i}, X_{j + 3, i},$ and $Y_i$ are functions on F$_{\infty}$.
We first note that all of the coefficients $ A, B, C, T_{i + 1, i}, X_{i + 2, i}, X_{j + 3, i}, Y_i$
must be functions on F$_{k + 3}$ in order for $\Phi$ to be closed, which in turn justifies the decomposition of $\bigwedge^2 \mathcal{I}_{\infty}$ in \emph{Theorem 2}. Moreover, the coefficient
of the highest weight $2 k + 3$ under the representation corresponding to the connection form $\beta$
of the original G structure F $\to$ M, T$_{k + 2, k + 1}$, turns out to indicate exactly when $\Phi \in \mathcal{H}_{k + 2} \subset \mathcal{H}_{k + 3}$. 

Throughout this section, $\Phi$ will always be assumed to be closed, and for the reasons
to follow, we call $A$, $B$, $C$, and $ \{ \, T_{i + 1, i} \, \}_{i = 0}^{\infty}$ the \emph{principal coefficients}.

\begin{theorem} For $\Phi \in \mathcal{H}_{k + 3}$ with $k \geq -1$,
 \begin{align} 
    &1. \; \; \Phi \in \mathcal{H}_{k + 2} \subset \mathcal{H}_{k + 3} \quad \mbox{iff} \quad T_{k + 2, k + 1} \, = 0. \notag \\ 
    &2. \, \, dT_{k + 2, k + 1} \,  \equiv T_{k + 2, k + 1} \,  ( \, \, 2\alpha - ( 2k + 3 ) \, \beta
                    \, + \, 2 N \, \theta_0 \, - 2 K \theta_1 \, )  \mod  \, \omega^2.     \notag
 \end{align}
\end{theorem}
\begin{corollary}
If d$( \alpha + N \, \theta_0 - K \, \theta_1 \, ) \, \, \ne 0$, every conservation law belongs to $\mathcal{H}_1$.
\end{corollary}
This corollary follows from the fact the differential system is time independent and d$\beta$ = 0. In local coordinates, the equation (1) satisfies
\begin{align}
   d( \alpha + N \, \theta_0 - K \, \theta_1 \, ) \, = \, 0 
\end{align}
iff $g(x, u, p, q)$ is at most quadratic in $q$, and $f(x, u, p)$ and $g(x, u, p, q)$ satisfy some fifth order partial differential relations. We also mention that if $K$ is not $0$, there is not any conservation laws in $\mathcal{H}_1$.

\vspace{1pc}
\noindent \emph{Proof of the Theorem}

\noindent We assume $k \geq 1$. The cases $k = -1,$  $0$ can be checked by direct computations.

1. From the results in the preceding section, we need to show
\begin{center}
   $ X_{k + 2, k} = \, Y_{k + 2} = \, X_{k + 3, i} = \, 0 \quad \mbox{for} \quad 0 \leq i \leq k $
\end{center}
assuming  $T_{k + 2, k + 1} = 0.$ Rather than getting into the detail of the computation, we present a brief sketch of the arguments involved. 

(a) $d\Phi \equiv 0 \quad \mod \quad \omega^2, \, \pi_0, \, ... \, \pi_k \, \, $ gives $X_{k + 2, k} = X_{k + 3, k} = 0.$

(b) $d\Phi \equiv 0 \quad \mod \quad \pi_0, \, ... \, \pi_k \, , \theta_{k + 1}, \, \theta_{k + 2} \,  \, $ gives $Y_{k + 2} = 0.$ Here we note that $d\theta_i \, \equiv \, d\pi_i \: \mod \: \omega^2.$

(c) First,
\begin{center}
$d\Phi \equiv 0 \quad \mod \quad \omega^2, \, \theta_0, \, ... \, \widetilde{\theta_k}, \,  \theta_{k + 1}, \, \sigma_{k + 2} \, \, $ 
\end{center}
gives $X_{k + 3, k - 1} = 0.$ 
Inductively,
\begin{center}
$d\Phi \equiv 0 \quad \mod \quad \omega^2, \, \theta_0, \, ...\, \widetilde{\theta_i}, \, \theta_{i + 1}, \, \theta_{i + 2}, \, \sigma_{i + 3}, \, ... \, \sigma_{k + 2} \, \, $ 
\end{center}
gives $X_{k + 3, i - 1} = 0 \quad \mbox{for} \quad k - 1 \geq i \geq 1.$ 
\noindent
Other than the structure equations in \emph{Theorem 1}, we make use of the fact
\begin{align}
d\sigma_i \, &\equiv \, dr_i \wedge \omega^1 \quad \mod \quad \alpha, \, \beta, \, \theta_0, \, 
                                                    \theta_1, \, \theta_2, \notag \\
  &\equiv \, 0  \quad \mod \quad \alpha, \, \beta, \, \theta_0, \, ... \, \theta_{i - 2}, 
                                      \, \sigma_{i - 1}, \notag
\end{align}
for $r_i$ is a function on F$_{i - 1}$ that is time independent.

2. We only show 
\begin{center}
$dT_{k + 2, k + 1} \equiv 0 \quad \mod \quad \alpha, \, \beta, \, \omega^2, \, \theta_0, \, \theta_1. $
\end{center}
The rest of the assertion follows easily from this. Put
\begin{align}
dT_{k + 2, k + 1} \, &\equiv \, T_{k + 2, k + 1}^{-1} \omega^1 \, + \, T_{k + 2, k + 1}^0 \theta_0 \, + . . . +  T_{k + 2, k + 1}^{k + 2} \theta_{k + 2} \, + \, T_{k + 2, k + 1}^{k + 3}  \sigma_{k + 3} \notag \\
  &  \mod  \alpha, \, \beta, \, \omega^2 \notag
\end{align}
to express the covariant derivative of $T_{k + 2, k + 1}$. 

First, $d\Phi \equiv 0 \quad \mod \quad \omega^2, \, \pi_0, \, ... \, \pi_k \, \, $ gives 
\begin{center}
$T_{k + 2, k + 1}^{k + 3} = 0, \, X_{k + 3, k} \, = \, - T_{k + 2, k + 1}, \: X_{k + 2, k} \, = \, - T_{k + 2, k + 1}^{-1}.$
\end{center}
Upon successively evaluating $d\Phi \, \equiv 0$ $\, \mod$ $\:$ $\pi_0,$ $\, ... \,$ $\pi_k,$ $\, \theta_{k + 1},$ $\, \theta_{k + 2}$ $\,$  and $ \, \mod \:$ $\pi_0,$ $\,$ $...$ $\,$ $\pi_k,$ $\, \theta_{k + 1}, \, \omega^1 \, $ with these relations, we get
\begin{align}
Y_{k + 2} \, = \, p_{k + 3} T_{k + 2, k + 1}^{-1}, \, X_{k + 3, k - 1} \, &= \, - T_{k + 2, k + 1}^{-1}. \notag
\end{align}
Now, $d\Phi \equiv 0 \quad \mod \quad \omega^2, \, \pi_0, \, ... \, \pi_{k - 1}, \, \theta_{k + 1}, \, \theta_{k + 2} \, $ gives
\begin{center}
$T_{k + 2, k + 1}^{-1} = 0$,
\end{center}
and we have
\begin{align}
dT_{k + 2, k + 1} \, &\equiv  \, T_{k + 2, k + 1}^0 \theta_0 \, + . . . + T_{k + 2, k + 1}^{k + 2} \theta_{k + 2} \: \mod \: \alpha, \, \beta, \, \omega^2. \notag
\end{align}
Finally
\begin{align}
d ( dT_{k + 2, k + 1} ) &\equiv 0 \quad \mod \quad \alpha, \, \beta, \, \omega^2, \, \theta_0, \, . . . \, \theta_i \notag
\end{align}
for $i$ $=$ $k + 2, k + 1, ... , 2 \, $ implies $\, T_{k + 2, k + 1}^i \, = \, 0. $  \begin{flushright} $\square$ \end{flushright}

Since the equation $d\Phi = 0$ is a linear differential system for the coefficients 
$ A,$ $B,$ $C,$ $T_{i + 1, i},$ $X_{i + 2, i},$ $X_{j + 3, i},$ and $Y_i$, we have the following corollary.
\begin{corollary}
For a conservation law $\Phi \in \mathcal{H}_{k + 3}$, all of $X_{i + 2, i}, X_{j + 3, i},$ and $Y_i$
can be expressed in terms of the principal coefficients $A$, $B$, $C$, $ \{ \, T_{i + 1, i} \, \}_{i = 0}^{k + 1}$ and their
successive covariant derivatives.
\end{corollary}
Thus, for each $k + 3, \, k \geq -2$, there is a linear differential system 
\begin{center} $ D_{k + 3} = D_{k + 3}$ [ $\, T_{k + 2, k + 1}, \, T_{k + 1, k}, \, . . . , \, T_{1, 0}, \, A, \, B, \, C \,$  ] 
\end{center}
 on $\mathcal{H}_{k + 3}$ whose solutions correspond to the conservation laws in $\mathcal{H}_{k + 3}$. 

The theorem also suggests and justifies the following definition.
\begin{definition}
Let $\Phi \in \mathcal{H}_{\infty}$ be a conservation law. We define
\[ weight \: \Phi = | \Phi | = \max_{T_{i + 1, i} \ne 0, i \geq 0 } 2 i + 1. \]
For a conservation law $\Phi \in \mathcal{H}_1$, we set $| \Phi | = -1. $ \\
Also, we put 
\begin{align}
     \mathsf{C}_0 &= \, 0 \notag \\              
   \mathsf{C}_{k + 3} &= \, \{ \, \Phi \in \mathcal{H}_{k + 3} \, \mid  \, d\Phi = 0 \, \} \notag \\
             &= \, \{ \, conservation \, \, laws \, \, of \, \, weight \, \, at \, most \, \, 2k + 3 \, \} \notag \\
   \mathcal{C}_{k + 3} &= \, \mathsf{C}_{k + 3} / \mathsf{C}_{k + 2} \notag \\ 
             &= \, \{ \, conservation \, \, laws \, \, of \, \, weight \, \, 2k + 3 \, \}. \notag  
 \end{align}                
\end{definition}
\noindent \begin{corollary}
 \[  \mathcal{C}_{k + 3} = \, \{ \, \Phi \in \mathsf{C}_{k + 3} \, \mid  \,  \, \, T_{k + 2, k + 1} \ne 0 \, \}, \]
and
\begin{align}
   \mathcal{C}_{\infty} &= \, \lim_{k \to \infty} \mathsf{C}_{k + 3} \notag \\
                        &= \, \bigcup_{k = -2}^{\infty} \mathcal{C}_{k + 3}. \notag \\
                        &\notag 
 \end{align}
\end{corollary}

The results in this section can be summarized as follows. Along with the prolongation tower $\mbox{F}_{k + 3}$
\begin{center}
 $\mbox{F}_{\infty} \to \mbox{F}_{k + 3} \to \mbox{F}_2 = \mbox{F}$, 
\end{center}
we attach the subspaces $\mathcal{H}_{k + 3}$
\begin{center}
 $\bigwedge^2 \mathcal{I}_{\infty} \supset \mathcal{H}_{\infty} \to \mathcal{H}_{k + 3} \to \mathcal{H}_2 \to \mathcal{H}_1,$ 
\end{center}
and a sequence of linear differential systems $D_{k + 3}$
 \[ \quad \: \: \: \: \, \qquad D_{\infty}  \to D_{k + 3} \to D_2 \to D_1 \]
whose solutions are the conservation laws of weight at most $2k + 3$. The sequence $\{D_{k + 3}\}_{k = -2}^{\infty}$ enjoys the property
\[ D_{k + 3} = D_{\infty} \mid_{T_{i + 1, i} = 0 \, \, \, \mbox{for} \, \, \, i \geq k + 2} .\] 
 
We finally mention that the weight of a conservation law roughly corresponds to the number of $x$ derivatives of $u$ involved in local coordinate expression, which agrees with the fact the orders of the conservation laws of $KdV$ equation jump by two.

\section{Examples}  

We compute the conservation laws of weight at most 3 for two classes of differential systems. A differential system will be called type ($\: n_{-1}$, $n_1$, $n_3 \:$) if it has $n_i$ conservation laws of weight $i$. 

\subsection{ Nonlinear differential systems with $N = 0$ that have at least one conservation law of weight -1, 1, and 3 respectively}

A complete classification can be given for this class of differential systems, where the nonlinearity assumption allows us to close up the linear differential system $D_3$.  First of all, a necessary condition to have a conservation law of weight $-1$ is $K$ $=$ $0$, and imposing the condition $N = 0$ and (5) of \emph{Corollary 2}, the structure equations in \emph{Proposition 1} become
\[
    d \begin{pmatrix}
           \omega^2\\ \theta_0\\ \omega^1\\ \theta_1\\  s_2
      \end{pmatrix}
     =  \begin{pmatrix}
             0 \\
          - \theta_1 \wedge \omega^1 \\
             0 \\ 
          - s_2 \wedge \omega^1 \\
            ( G \theta_1 + L \theta_0 ) \wedge \omega^1     
        \end{pmatrix}
\]
on a section $\, \alpha \, = \, \beta \, = \, 0 \,$, since $\, d\alpha \, = \, d\beta \, = \, 0 \,$. It's easily verified that these equations imply the following simple local normal form
\begin{align}
    \omega^2 \, &= \, dt, \, \,  \omega^1 \, = \, dx \notag \\
    \theta_0 \, &= \, du - p dx \notag \\
     \theta_1 \, &= \, dp - q dx \notag \\
     s_2     \, &= \, dq + g(x,u,p) dx \notag
\end{align}
that corresponds to the evolutionary pde
\begin{equation} u_t \, = \, u_{xxx} \, + \, g(x, u, p). \notag \end{equation}
Notice the admissible point transformations at this stage are, up to scaling by constants, 
\begin{align}
   \bar{x} \, &= \, x \, + \, x_0 \notag \\
   \bar{u} \, &= \, u \, + \, \phi(x), \notag
\end{align}
where $x_0$ is a constant and $\phi(x)$ is an arbitrary function of $x$. \

Following the general theory developed earlier, a conservation law $\Phi$ $\, \in$ $\, \mathcal{H}_3$
takes the form
\begin{align}
   \Phi \, &= \, A ( \, \theta_0 \wedge \omega^1 \, + \, s_2 \wedge \omega^2 \, ) \, - \, B \theta_1
                          \wedge \omega^2 \, + \, C \theta_0 \wedge \omega^2  \notag \\
           & \, + T \, \theta_1 \wedge \theta_0 \notag \\
           & \, + P \, ( \, \theta_2 \wedge \theta_1 \, + \, s_3 \wedge \pi_0 \, ) , \notag
\end{align}
where 
\begin{align}
    \pi_0    \, &= \, \theta_0 \, \: \qquad \quad \:  - \, p_3 \, \omega^2 \notag \\
    \theta_2 \, &= \,  s_2 \, - \, p_3 \, \omega^1 \notag \\
    s_3      \, &= \, dp_3 . \notag
\end{align}
Then, a direct computation shows $\, d\Phi \, = \, 0 \,$ if 
\begin{align}
  dP \, &\equiv \, 0 \quad \mod \quad \omega^2  \notag \\
  T^3 \, &= \, T^2 \, = \, 0, \, \, T^{-1} \, = \, A^1 \, + \, g_u P\notag \\  
  A^3 \, &= \, 0, \: A^2 \, = \, T, \: A^{-1} \, = \, B   \notag    \\
  B^3 \, &= \, 0, \: B^2 \, = \, P^{-2} \, - \, A^1, \: B^{-1} \, = \, C \, - \, g_p A \, - \, p_3 P^{-2}  \notag \\
  C^3 \, &= \, 0, \: C^2 \, = \, A^0, \: C^1 \, = \, T^{-2} \, - \, B^0, \: C^{-1} \, = \,  g_u A \, + \, A^{-2}  \notag , 
\end{align}
where we express the covariant derivative of a coefficient $X$ by 
\begin{equation} dX \, = \, X^{-1} \omega^1 \, + \, X^{-2} \omega^2 \, + \, X^0 \theta_0 \, + \, X^1 \theta_1 \, + \, X^2 s_2 \, + \, X^3 s_3.  \notag \end{equation}
Note by \emph{Theorem 3}, 
\begin{align}
    \Phi \, &\in \, \mathcal{H}_2 \quad \mbox{if} \quad  P \, = \, 0  \notag \\
    \Phi \, &\in \, \mathcal{H}_{1} \quad \mbox{if} \quad  P \, = \, T \, = \, 0.  \notag
\end{align}

From now on, we give only the sketch of the procedure of overdetermined pde machinery that is applied to find the compatibility conditions to admit the desired conservation laws. When $X$, $Y$, . . . are the coefficients of the conservation laws and their covariant derivatives, the notation $ \{ \, X, \, Y, . . . \, \}$ for simplicity would mean an expression that is linear in $X,$ $Y,$  $...$  with the coefficients in $g(\, x, \, u, \, p \,)$ and its derivatives. For example, $g_u^2 P \, + \, g_p A^{-2}$ will be denoted by $\{ \, P, \, A^{-2} \, \}. $ We also mention if there exists a compatibility condition of the form
\begin{equation} \{ \, A, \, T, \, P \, \} \, = \, 0, \notag \end{equation}
all of the coefficients of $A,$ $\, T,$ and $P$ must be $0$ if the differential 
system has all the conservation laws of weight $-1$, $1$, and $3$.

\vspace{1pc}
\noindent
1. $d( d ( A ) ) \, \equiv \, 0 \: \mod \: \omega^2, \, \theta_0, \, \theta_1$ gives $ A^1 \, = \, \{ \, P, \, P^{-2} \, \}. $ 

   $d( d ( T ) ) \, \equiv \, 0 \: \mod \: \omega^2, \, \theta_0, \, \theta_1$ gives $ T^1 \, = \, 0 . $ 

   $d( d ( T ) ) \, \equiv \, 0 \: \mod \: \omega^2, \, \theta_0$ gives
$ T^0 \, = \, \{ \, P \, \} . $

   $d( d ( T ) ) \, \equiv \, 0 \: \mod \: \omega^1, \, \omega^2$ gives
$   \, g_{p p u} \, = \, 0, \, $
which implies
\begin{equation} g( x, u, p ) \, = \, h( x, u ) p \, + \, a( x, u ) \, + \, c( x, p ) \notag \end{equation}
for some functions $h, \, a,$ and $c$. 

But, $d( d ( T ) ) \, \equiv \, 0 \: \mod \: \omega^2$ gives
$ \,  a_{u u} \, = \, h_{x u}, \, $
which suggests writing
\begin{equation} g( x, u, p ) \, = \,  f( x, u )_x \, + \, c( x, p ). \notag \end{equation}

\noindent
2.  $d( d ( B ) ) \, \equiv \, 0 \: \mod \: \omega^2, \, \theta_0, \, \theta_1$ gives
$ B^1 \, = \, \{ \, P, \, T, \, A^0 \, \}. $ 

    $d( d ( A ) ) \, \equiv \, 0 \: \mod \: \omega^2$ gives
$ A^{0 -1} \, = \, \{ \, B^0, \, T \, \}, \: A^{0 1} \, = \, \{ \, P \, \}, \: A^{0 2} \, = \, \{ \, P \, \}. $

    $d( d ( B ) ) \, \equiv \, 0 \: \mod \: \omega^2, \, \theta_0, $ gives
$ B^0 \, = \, \{ \, P, \, P^{-2}, \, T, \, T^{-2}, \, A \, \}.  $

    $d( d ( B ) ) \, \wedge \, \omega^1 \, \wedge \omega^2 \, + \, \frac{1}{3} \, d ( d ( T ) ) \, \wedge \, \omega^1 \, \wedge \theta_0 \, = \, 0 $ gives
$  A^{0 0 } \, = \, \{ \, P, \, P^{-2}, \, T, \, A \, \}.  $

    $d( d ( A^0 ) ) \, \equiv \, 0 \: \mod \: \omega^1, \, \omega^2, \, P $ gives
$ c_{p p p p} \, = \, 0 \, = \, c_{p p p x}, $
which implies $c( x, p )$ is at most cubic in $p$. 

Let's consider the case $c( x, p )$ is quadratic in $p$, for example, and write
\begin{equation} g ( x, u, p ) \, = \,   f( x, u )_x \, + \, c_2( x ) p^2 / 2 \, + \, c_1( x ) p  \notag \end{equation} 
after adding some function of $x$ to $u$ if necessary.

\vspace{1pc}
\noindent
3.  $d( d ( A^0 ) ) \, \equiv \, 0 \: \mod \: \omega^2, \, \theta_1, \, P, \, T $ gives
$ c_2( x ) \, = \, 0. $

    $d( d ( A^0 ) ) \, \equiv \, 0 \: \mod \: \omega^1, \, \omega^2, $ gives
$ f_{u u u u} \, = \, 0 \, = \, f_{u u u x}, \: $
and $f( x, u )$ is at most cubic in $u$. 

Let's again consider the case $f( x, u )$ is quadratic in $u$ and write
\begin{equation} g( x, u, p ) \, = \, ( \, f_2( x ) \, u^2/2 \, + \, f_1( x ) \, u \, + f_0( x ) \, )_x \, + \, c_1( x ) \, p , 
\notag \end{equation} 
assuming $\, f_2( x ) \, \ne \, 0 \,$ for nonlinearity. Again by adding a suitable function of $x$ to $u$, we may write
\begin{equation} g( x, u, p ) \, = \, ( \, f_2( x ) \, u^2/2  \, + f_0( x ) \, )_x \, + \, c_1( x ) \, p . \notag \end{equation}

\noindent
4.  $d( d ( A^0 ) ) \, \wedge \, \omega^2 \, + \, \frac{1}{3} \, d ( d ( T ) ) \, \wedge \, \omega^1  \, = \, 0 $ gives 
$\,   P^{-2} \, = \, \{ \, P \, \}, \, f_2' \, = \, 0. $

Hence $f_2$ is a nonzero constant and after scaling, we may write
\begin{equation} g( x, u, p ) \, = \, ( \, u^2/2 + f_0( x ) \, )_x \, + \, c_1( x ) \, p . \notag \end{equation} 
   
   $ d( d( P ) ) \, = \, 0 $ gives
$\,  c_1'' \, = \, 0, \,$
and we put $ \, c_1( x ) \, = \, r_1 x \, + \, r_0 \, $ with constants $r_1$, $r_0$. 

   $d( d ( B ) ) \,  \wedge \omega^2 \, + \, \frac{1}{3} \, d ( d ( T ) ) \, \wedge \, \theta_0  \, = \, 0 $ gives
$\,  C^0 \, = \, \{ \, P, \, A, \, A^0 \, \}. \, $

   $d( d( C ) ) \, \wedge \, \omega^2 \, + \, d( d ( A ) ) \, \wedge \, \omega^1 \, + \, \frac{2}{3} \, d ( d ( T ) ) \, \wedge \, \theta_1  \, = \, 0 $ gives
$\,  A^{0 -2} \, = \, \{ \, P, \, T, \, T^{-2}, \, A^0, \, B \, \}. \,$ 

   $d( d ( A^0 ) ) \, \wedge \, \omega^1 \, - \, \frac{1}{3} \, ( u \, + \, c_1( x ) ) \, d ( d ( T ) ) \, \wedge \, \omega^1  \, = \, 0 $ gives
$\,  T^{-2} \, = \, \{ \, P, \, T \, \}. \,$

Finally, $d( d ( T ) ) \, = \, 0 $ gives
$\,  f_0'' \, = \, - 2 r_1^2 \,$
and we put
\begin{equation} g( x, u, p ) \, = \, ( \,  u^2/2  \, + ( \, - r_1^2 \, x^2 \, + \, b_1 \, x  \, ) )_x \, + \, ( \,  r_1 x \, + \, r_0 \, )\, p  \notag \end{equation} 
for some constant $b_1$.

\vspace{1pc}
\noindent
5. $ d( d ( A^0 ) ) \, = \, 0\, $ gives  $\, C \, = \, \{ \, P, \, T, \, A, \, A^0, \, B \, \}, $
and the exterior derivative of this relation gives
$A^{-2} \, = \, \{ \,  P, \, T, \, A, \, A^0, \, B \, \}.$ Now, by successively taking $d ( d ( A ) ), \, d ( d ( B ) )$, the remaining variables $B^{-2}$ and $C^{-2}$ are determined in terms of the free variables  $P$, $T$, $A$, $B$, $A^0$ , and the linear differential system for these free variables becomes completely integrable. Since
\begin{align}
dP \, &\equiv \, 0 \quad \mod \quad P \notag \\
dT \, &\equiv \, 0 \quad \mod \quad P, \, T   , \notag
\end{align}
this class of differential systems (equations) are type $ ( \, 3, \, 1, \, 1 \, )$. 

\vspace{1pc}
The computations similar to the one demonstrated above determine the following list of (nonlinear) differential systems (equations) and their types.
\begin{align}
  u_t \, &=  u_{xxx}  +  ( \frac{1}{2}u^2  -  r_2^2  x^2  +  r_1  x  )_x  +  r_2 x  p                                    \qquad ( 3, 1, 1 )  \\
  u_t \, &=  u_{xxx}  \pm ( \frac{1}{6}u^3 + r_1 x u + r_0 u)_x \qquad \quad  \qquad ( 2, 2, 1 ) \notag \\
          & \notag \\
  u_t \, &= \, u_{xxx}  +  ( m' \cos(u) + n' \sin(u) )_x  + \frac{1}{8}p^3  + b p  \qquad \; \; \; ( 2, 1, 2 ) \notag \\
  u_t \, &= \, u_{xxx}  +  ( m \cosh(u) + n \sinh(u) )_x  - \frac{1}{8}p^3  + b p \qquad \; ( 2, 1, 2 ) \notag \\
          & \notag \\
  u_t \, &= \, u_{xxx}  \pm  ( \exp(u) )_x - \frac{1}{8}p^3 + ( r_1 x + r_0 ) p + 2 r_1 \qquad ( 2, 1, 2 ) \notag \\
  u_t \, &= \, u_{xxx} \pm  \frac{1}{6}p^3+ ( r_1 x + r_0 ) p + b \qquad \qquad \qquad \qquad \; ( 2, 1, 2 ) \notag 
\end{align}
where $r_2, r_1, r_0, b, m', n', m, n$ are constants such that $m^2$ $\ne$ $n^2$.

\subsection{ Flow of the curve in the plane by the derivative of its Finsler curvature with respect to the arc length  }  

Consider the two dimensional abelian group $R^2$ with a translation invariant Finsler structure, i.e., a translation invariant norm on the tangent vectors.[B] \; It is known there exists a coframe $ \{ \: \eta^1, \, \eta^2, \, \eta^3 \: \} $ on the unit circle bundle $\mathcal{B} \, \to \, R^2$ with the structure equations
\begin{align}
   d\eta^1 \, &= \, \eta^3 \wedge \eta^2   \notag  \\
   d\eta^2 \, &= \, - \eta^3 \wedge ( \, \eta^1 \, - \, I \eta^2 \, )   \notag  \\
   d\eta^3 \, &= \, 0,   \notag 
\end{align}
where $I$ is the fundamental invariant satisfying 
\begin{equation} dI \, \equiv \, 0 \quad \mod \quad \eta^3. \notag \end{equation} 
Here $\eta^1$ is the Hilbert form and $\eta^3$ is the canonical pseudo connection form. 

Given an immersed curve $\gamma$ in the plane, $\bar{\gamma}^{*} ( \, \eta^1 \, ) \ne 0 \, $, we define as usual the curvature of $\gamma$ by
\[  \bar{\gamma}^{*} ( \, \eta^3 \, )  \, = \, k \bar{\gamma}^{*} ( \, \eta^1 \, ), \]
where $\bar{\gamma}$ is the canonical lift of $\gamma$ to $\mathcal{B}$. We also define $k_1$, the derivative of $k$ with respect to the arclength, by
\[ dk \, = \, k_1 \, \bar{\gamma}^{*} ( \, \eta^1 \, ). \] 

The \emph{ $k_1$ flow } we are interested in is the following. Given $\gamma$ in $R^2$, first lift it to $\mathcal{B}$, $\bar{\gamma}$. Then flow $\bar{\gamma}$ along the Finsler normal, or the vector dual to $\eta^2$, in the amount of $k_1$. The image of this flow under the projection $\mathcal{B} \, \to \, R^2$ will be called $k_1$ flow of $\gamma$.

This $k_1$ flow can be expressed in terms of a differential system on M$^5 \, = \, $ $ \, \{ \, t \, \}$ $\times$ $\{ \, k \, \}$ $\times$  $\mathcal{B}$ as follows.
\[ \mathcal{I} \, = \, \{ \, \theta_0 \wedge \omega^1 \, + \, \sigma_2 \wedge \omega^2, \, \theta_1 \wedge \omega^2, \, \theta_0 \wedge \omega^2 \, \} \cup \{ \, \theta_1 \wedge \theta_0 \, \}, \]
where
\begin{align}
   \omega^2 \, &= \, dt \notag \\
   \theta_0 \, &= \, \eta^2 \notag \\
   \omega^1 \, &= \, \eta^1 - I \eta^2 \notag \\
   \theta_1 \, &= \, ( \, \eta^3 - k \eta^1 \, ) \, - \, k I \eta^2 \notag \\
   \sigma_2 \, &= \, dk,    \notag
\end{align}
with the independence condition $\omega^1 \wedge \omega^2 \ne 0$. It is clear the integral manifold of $\mathcal{I}$ on which $\omega^1 \wedge \omega^2 \ne 0$ is the graph over the canonical lift of a $k_1$ flow. 

The computation similar to that of \emph{Example (1)}, i.e., repeated application of the identity $d^2 \, = \, 0$, yields the following results.

\vspace{1pc}
(a) case $\, I \, = \, 0$ \\
\noindent In this case, the Finsler structure is the flat Riemannian metric on $R^2$. The differential system is type $( \, 4, \, 1, \, 1 \, )$. If the initial curve $\gamma$ is simply closed, the first five conservation laws represent the preservation under the flow of 
\[ \int_D c_1 + c_2 x + c_3 y + c_4 ( x^2 + y^2 ) \: dx \wedge dy, \]
and the length of the curve, where ($x, y$) are the orthonormal coordinates of the plane, $D$ is the region enclosed by the curve, and $c_i'$s are constants. We also mention $\gamma$ admits a self similar $k_1$ flow if 
\[ k_1^2 + \frac{1}{4} k^4 + a_2 k^2 + a_1 k + a_0 = 0  \]
for some constants $a_i$. It is interesting to note that upon imposing $a_2$ $=$ $a_1$ $=$ $0$,
the above ordinary differential equation for $k$ becomes the \emph{Euler - Lagrange} equation for the functional
\[ \int_{\gamma} \, \frac{1}{2} k^2 \, . \]
In [ChT], it is shown the curvature of the $k_1$ flow satisfies $mKdV$ equation under a suitable local coordinates system. 

\vspace{1pc}
In case $\, I \, \ne \, 0$, the differential system does not have any conservation laws of weight -1. Below are the Finsler structures whose $k_1$ flow have at least one conservation law of weight 1 and 3.

\vspace{1pc}
(b) case $\, I \, \ne \, 0, \, dI \, = \, ( \, - \frac{3}{2} \, + \, \frac{1}{3} \, I^2 \, ) \, \eta^3$ 

\noindent It is type $( \, 0, \, 3, \, 2  \, )$ and also has one conservation law of weight 5. 
There does not exist a Finsler structure with this structure equation.[B]

\vspace{1pc}
(c) case $\, I \, \ne \, 0, \, dI \, = \, I_3 \eta^3 \,\ne \, ( \, - \frac{3}{2} \, + \, \frac{1}{3} \, I^2 \, ) \, \eta^3, \: dI_3 \, = \, ( \, - I + \frac{2}{9} \, I^3 \, ) \, \eta^3$ 

\noindent It is type $( \, 0, \, 3, \, 1 \, )$ and also has one conservation law of weight 5.
\[ \]

We finally mention that each of the differential systems (equations) presented in the list (6) of \emph{Example (1)} 
and  (a), (b), (c) of \emph{Example (2)}  also has exactly one conservation law of weight 5. In fact, 
after the computation of the examples, it is tempting to conjecture the following. 

\vspace{1pc}
\emph{ Consider a differential system that locally corresponds to the nonlinearizable evolutionary pde of the form
\[ u_t \, = \, f( \, x, \, u, \, u_x \,) \, u_{xxx} \, + \, g( \, x, \, u, \, u_x, \, u_{xx} \,).  \]
If the differential system has the conservation laws of weight $ \, -1, \, 1, \, $ and $3$, then it has exactly one 
conservation law of weight $\, 2 k + 3 \,$ for each $k \, \geq \, 1$. In general, if the differential system has 
two conservation laws of distinct weights $\geq \, 1$, 
then it has an infinite sequence of conservation laws of distinct weights. } 

\vspace{2pc}
\noindent \textbf{\Large{References }} 

\vspace{1pc}
\noindent
[B] R. Bryant, \emph{Finsler surfaces with prescribed curvature conditions}, preprint (1995).

\vspace{1pc}
\noindent
[BG1] R. Bryant and P. A. Griffiths, \emph{Characteristic cohomology of differential systems I},  J. Amer. Math. Soc. \textbf{8} (1995), 507-596. 

\vspace{1pc}
\noindent
[BG2] \underline{\quad }\underline{ }, \emph{Characteristic cohomology of differential systems II}, Duke Math. J. \textbf{78} (1995), 531-676. 

\vspace{1pc}
\noindent
[BG3] \underline{\quad }\underline{ }, \emph{Reduction for constrained variational problems and $\int \frac{1}{2}\kappa^2 ds$}, Amer. J. Math. \textbf{108} (1986), 525-570. 

\vspace{1pc}
\noindent
[ChT] S. S. Chern and K. Tenenblat, \emph{Foliations on a surface of constant curvature and the modified Korteweg - De Vries equations}, J. Diff. Geom. \textbf{16} (1981), 347-349.

\vspace{1pc}
\noindent
[F] K. Foltinek, \emph{Quasilinear third order scalar  evolution equations and their conservation laws}, Thesis, Duke university (1996). 

\vspace{1pc}
\noindent
[K] O. V. Kaptsov, \emph{Classification of evolution equations by conservation laws}, Functional Analysis and its Appl. \textbf{16} no.1 (1982), 59-61. 

\vspace{1pc}
\noindent
[M] A. V. Mikhailov, A. B. Shabat, and V. V. Sokolov, \emph{The symmetry approach to classification of integrable equations}, What is integrability ?, Springer Series in Nonlinear Dynamics, Springer-Verlag, Berlin (1991).

\vspace{2pc}
\noindent Duke University, \\
Durham, NC 27708 \\
\emph{E-mail address} : \textbf{ship@math.duke.edu}

\end{document}